\documentclass[journal]{IEEEtran}

\usepackage{cite}
\usepackage{color}

\ifCLASSINFOpdf
  \usepackage[pdftex]{graphicx}
  \graphicspath{{../pdf/}{../jpeg/}}
\else
\fi

\usepackage{amsmath}
\usepackage{amsfonts}
\usepackage{amssymb}
\usepackage{cases}
\newtheorem{lemma}{Lemma}

\newtheorem{assumption}{Assumption}
\newtheorem{property}{Property}
\newtheorem{remark}{Remark}

\usepackage{multirow}
\usepackage{threeparttable}
\usepackage{circuitikz}
\ctikzset{bipoles/length=1cm}

\usepackage{url}

\allowdisplaybreaks[4]

\usepackage{enumerate}
\usepackage{color}

\newcommand{\lj}[1]{\textcolor{black}{{}#1}}

\usepackage{comment}
\usepackage{ifthen}

\newboolean{showcomments}
\setboolean{showcomments}{true}
\newcommand{\fliu}[1]{\ifthenelse{\boolean{showcomments}}
	{ \textcolor{red}{(FL:  #1)}}{}}
\newcommand{\jli}[1]{\ifthenelse{\boolean{showcomments}}
	{ \textcolor{blue}{(JL:  #1)}}{}}

\begin{document}

\title{Resilience Control of DC Shipboard Power Systems}

\author{Jia~Li,
        Feng~Liu,~\IEEEmembership{Member,~IEEE,}
        Ying~Chen,~\IEEEmembership{Member,~IEEE,}
        Chengcheng~Shao,~\IEEEmembership{Member,~IEEE,}
        Guanqun~Wang,~\IEEEmembership{Member,~IEEE,}
        Yunhe~Hou,~\IEEEmembership{Senior Member,~IEEE}
        Shengwei~Mei,~\IEEEmembership{Fellow,~IEEE}
\thanks{
\lj{This work was supported in part by the National Natural Science Foundation of China (51677100, 51621065), in part by the major basic research in the National Security of China (613294). \emph{(Corresponding Author: Feng Liu)}}}
\thanks{J. Li, F. Liu, Y. Chen and S. Mei are with the State Key Laboratory of Power Systems, Department of Electrical Engineering, Tsinghua University, Beijing 100084, China.} 
\thanks{C. Shao is with School of Electrical Engineering, State Key Laboratory on Electrical Insulation and Power Equipment, Shanxi Province Key Laboratory on Smart Grid in Xi'an Jiaotong University, Xi'an, 710049, China.}
\thanks{G. Wang is with Burns \& McDonnell, Houston, TX 77027 USA.}
\thanks{Y. Hou is with the Department of Electrical and Electronic Engineering, University	of Hong Kong, Hong Kong, China.}
        }

\maketitle

\begin{abstract}
Direct current (DC) network has been recognized as a promising technique, especially for shipboard power systems (SPSs). Fast resilience control is required for an SPS to survive after faults. Towards this end, this paper proposes the indices of survivability and functionality, based on which a two-phase resilience control method is derived. The on/off status of loads are determined in the first phase to maximize survivability, while the functionality of supplying loads are maximized in the second phase. Based on a comprehensive model of a DC shipboard power systems (DC-SPS), the two-phase method renders two mixed-integer non-convex problems. To make the problems tractable, we develop second-order-cone-based convex relaxations, thus converting the problems into mixed-integer convex problems. Though this approach does not necessarily guarantee feasible, hence global, solutions to the original non-convex formulations, we provide additional mild assumptions,
which ensures that the convex relaxations are exact \lj{when line constraints are not binding.}
In the case of inexactness, we provide a simple heuristic approach to ensure feasible solutions. Numerical tests empirically confirm the efficacy of the proposed method.
\end{abstract}


\IEEEpeerreviewmaketitle

\section{Introduction}

\IEEEPARstart{D}{riven} by the increasingly strict government regulation of emissions and customers' ever-increasing fuel-efficiency requirements, all-electric ship (AES) turns out to be an emerging trend, which exploits an electrical propulsion system instead of the conventional mechanical system. In such a circumstance, a direct current (DC) power system has several inherent advantages over an AC system, such as \cite{Jin2016a}:
\begin{itemize}
  \item Replacing bulky ferromagnetic transformers with compact power electronic converters;
  \item Easier implementation of parallel connection or disconnection for DC power sources;
  \item Avoiding synchronization problems;
  \item Exempting from reactive power flows;
  \item Eliminating harmonic and imbalance problems.
\end{itemize}

Additionally, considering the specific needs of shipboard power systems (SPSs), DC networks could bring a broad range of advantages \cite{Jin2016a,Reed2012}. A DC network eliminates bulky low-frequency transformers and reduces the rating of switchgear, thus reducing the occupied space and overall weight of the whole system. Moreover, a DC network could provide stronger survivability, better limitation of fault current, and higher reconfiguration capability. In addition, the integration of advanced high-speed, high-efficiency diesel generation (i.e., gas turbine generation) could also be easily achieved within a DC network, which could effectively improve the fuel efficiency of system operation. Hence, DC shipboard power systems (DC-SPSs) keep gaining an increase in research interests. 

Being entirely self-contained, SPSs may suffer severe consequences from faults. After faults, fast resilience control for service restoration can greatly facilitate the survivability of SPSs. In the literature, different network reconfiguration approaches have been proposed for resilience control in SPSs. A modified fixed-charge network flow method is proposed in \cite{Butler1999} , and further improved to consider more constraints \cite{Butler2001} and incorporate geographical fault information \cite{Butler-Purry2004}. Distributed approaches, such as multi-agent system technique is applied to solve SPS reconfiguration problems, which relies only on local information \cite{Butler-Purry2005,Solanki2005,Huang2007,Feliachi2006a}. Moreover, intelligent algorithms, such as genetic algorithm \cite{Jing2009,Shariatzadeh2016}, particle swarm optimization \cite{Kumar2007,Mitra2011,Shariatzadeh2016}, ant colony optimization \cite{Hari2010} are also employed. However, these works focus on AC SPSs. 

With respect to DC SPSs, an optimal method for restoring power to loads, including islanding, is proposed in \cite{Khushalani2008}. The service restoration problem is formulated for both balanced and unbalanced SPS. A dynamic formulation and a genetic-algorithm-based static implementation of a damage control method at the DC zonal SPS are presented in \cite{Amba2009}. In \cite{Bose2010,Bose2012}, a convex form of the optimal reconfiguration problem is obtained by an affine transformation of the nonlinear equality constraints using Newton's power flow method.
Cumulative distribution function of the power delivered to loads is presented to showcase the system robustness against random fault scenarios. Dynamic technique is employed for automatic reconfiguration in \cite{Das2013}, which produces not only the final reconfiguration, but also the correct order in which the switches are to be changed.

Quadratic relaxations of power flow equations for radial and meshed networks are extended to include topology changes and \lj{formulate} reconfiguration problems \cite{Hijazi2015}. 
In \cite{Gan2014}, a second-order cone program (SOCP) was first suggested to exactly relax the OPF problem in \emph{grid-connected} DC distribution networks, with strict proofs and comprehensive discussions. It laid \lj{the} foundation for solving \lj{optimal} power flow (OPF) problems in DC distribution networks, and was further extended in \cite{Li} to a broader class of DC networks, including \emph{stand-alone} DC microgrids. However, discrete operations, such as line switching and re-sectionalization, have not yet been considered in both works.
In this context, this paper proposes a two-phase resilience control method for DC SPSs. To this end, indices of survivability and functionality are first proposed. The first phase is to maximize the survivability, where the network is reconfigured, probably with a certain amount of load shedding. The second phase is to maximize the functionality by adjusting the supplied loads. To formulate DC SPSs with line switching actions, spanning tree constraints are derived. The action of bus-tie breakers is modeled so as to exploit the flexibility of re-sectionalization. Additionally, both AC-DC and DC-DC converters are modeled from the view of DC network side. Since the established optimization problems are inherently mixed-integer and non-convex, it is difficult to guarantee the global optimality. Therein, convex relaxation is \lj{applied} to convert them into mixed-integer second-order cone programming (MISOCP) problems. Under mild assumptions, a convex relaxation method is proposed. When the current limit constraint is inactive, the relaxation is proven to be exact. In this context, the global optimal solution is guaranteed theoretically. However, when the current limit constraint is active, the relaxation may not be exact. Thus, a heuristic method is devised to give an approximate solution if the original solution is not exact.

It should be noted that, whereas our work is focused on shipboard power systems, other bulk DC power systems or DC microgrids can also benefit from our methodology, including:
\begin{enumerate}
\item \textbf{Formulation:} The line switching options are modeled, enabling us to exploit the flexibility of network topology via line switching and re-sectionalization. This formulation can also be applied to other bulk DC power systems. The spanning tree constraints can be used to model other DC power networks with similar topologies, or they may be simply removed to model meshed topologies.

\item \textbf{Solution approach:} The proposed MISOCP-based convex relaxation is independent of topology and operating mode (grid-connected or stand-alone) of the power network. Thus, it can also be applied to various kinds of bulk DC power systems considering topology changes.
\end{enumerate}

The rest of this paper is organized as follows. The indices of survivability and functionality are proposed in Section \ref{sec:index}. The two-phase resilience control method is presented in Section \ref{sec:formulation}. The MISOCP solution approach is given in Section \ref{sec:solution}. Numerical studies are provided in Section \ref{sec:case}. Section \ref{sec:conclusion} concludes the paper with remarks.

\section{Indices and Framework}\label{sec:index}


In the resilience control of SPSs, the first priority is to retain the critical loads in order to survive. Then these critical loads are supplied partly or fully,  according to the needs of different functions. 
Note that the two objectives may not be consistent usually. In line with the decision-making process, we propose two key indices, namely survivability and functionality. Then we derive a two-phase framework of resilience control for DC-SPSs. In Phase I, we maximize the survivability index by adjusting the switches of \lj{lines and} loads, so that the critical loads are guaranteed being connected to the SPS in order of priority. In Phase II, with the status of load switches obtained in Phase I, we further maximize the functionality index  by adjusting the amount of loads and line switches. 

Mathematically, a power network can be abstractly depicted by a graph, denoted by $\mathcal{G}:=(\mathcal{N},\mathcal{E})$. Here $\mathcal{N}:=\{1,\cdots,n\}$ stands for the set of all buses and $\mathcal{E}$ denotes the set of all lines in the network. 
For each bus $i\in \mathcal{N}$, denote $p_i$ as its power injection (either generation or load). Let $\underline{p}_i$ and $\overline{p}_i$ indicate the lower and upper bound on $p_i$, respectively.
Let $\mathcal{N}^D$ denote the set of load buses, and $\mathcal{N}^G$ the set of generator buses.

\subsection{Survivability Index}
In this subsection, a survivability index is proposed to guarantee that most important loads are restored priorly. 

In traditional methods, weighting factors are introduced into the objective function to indicate the priority of loads in reconfiguration. A simple and effective method is proposed in \cite{Nagarajan2016,Barnes2017}, which uses some constraints to guarantee that a minimum predefined fraction of loads is satisfied. In this paper, we assume the minimum fraction of loads that must be served is unknown before a fault occurs, since the fault may be so severe that the predefined requirement can not be satisfied. Thus, an alternative method is proposed to circumvent such a problem.

Additional binary variables are introduced here, as explained. For each load, $\eta_i$ denotes a binary variable indicating whether load $i (i\in \mathcal{N}^D)$ is switched on ($\eta_i=1$) or  off ($\eta_i=0$). 
Then the survivability index can be built on the rule that the loads are satisfied in order of priority, regardless of load capacity. 
Let $\mathcal{M}:=\{1,\cdots,M\}$ denote the set of priority levels. 
For each priority level $m\in\mathcal{M}$, $\phi_m$ denotes the number of loads pertaining to it. Each load is assigned with a weighting factor $\kappa_i^m$.  Assume that the smaller number of $m$ indicates the higher priority. Then the normalized survivability index is defined as follows. 
\begin{equation}\label{eq:omegas}
\Omega_{sur} := \frac{\sum_{i\in\mathcal{N}^D}\kappa_i^m\eta_i}{\sum_{i\in\mathcal{N}^D}\kappa_i^m}
\end{equation}
where the weighting factor of each priority level satisfies
\begin{equation}\label{eq:weight}
\kappa_i^m > \sum_{r=m+1}^M \phi_r\kappa_i^r.
\end{equation}
The denominator in \eqref{eq:omegas} is a constant for normalization so that $\Omega_{sur}\in [0,1]$. 
Inequality \eqref{eq:weight} implies that the weighting factor of a higher priority $m$ must be greater than the summation of all the weighting factors associated with lower priorities, i.e., $k+1,\cdots,M$. Therefore, when any load with a higher priority is switched on, $\Omega_{sur}$ will be greater than that when all the loads with  lower priorities are switched on.

\subsection{Functionality Index}

When all the loads in a system can be satisfied at their maximum power, it is known as the function of the system is fully fulfilled. To quantify the (weighted) proportion of load that can be supplied, an additional index, which is referred to as functionality index, is defined as follows.
\begin{equation}\label{eq:omegaf}
\Omega_{fun} := \frac{\sum_{i\in\mathcal{N}^D}\eta_i \lambda_i p_i}{\sum_{i\in\mathcal{N}^D}\eta_i \lambda_i \underline{p}_i}
\end{equation}
where $\lambda_i$ are weighting factors in terms of different function requirements. Binary variables $\eta_i$ ensure that only the loads switched on are taken into account in the index. The denominator in \eqref{eq:omegaf} is a constant for normalization so that $\Omega_{fun}\in [0,1]$.

\subsection{A Two-phase Framework of Resilience Control}

Leveraging the two indices above, a two-phase  framework of resilience control for DC-SPSs can be formed as below.

\subsubsection{\textbf{Phase I}}
In this stage, the survivability index, $\Omega_{sur}$, is maximized by adjusting load and line switches, so that the loads can be switched on strictly subject to their priorities.
\subsubsection{\textbf{Phase II}} 
With the status of load switches obtained in Phase I, the functionality index, $\Omega_{fun}$, is further maximized by adjusting the amount of loads and line switches, so that the system functions can be fulfilled as much as possible.

\begin{remark}
	In this paper, we devise a two-phase framework, instead of combining the two phases into a whole problem, due to the following reason. In the resilience control design of SPSs, the two objectives, i.e., survivability and functionality, are not equally important in practice. The survivability has the \emph{highest} priority which must be satisfied first, while the functionality has the second which should not influence the survivability. However, if the two phases are forcibly merged into one problem, the functionality index has to be enforced upon the objective function of the Phase-I problem. This treatment will undermine the optimal solution of the Phase-I problem, making it deviate from the global optimum.  On the other hand, mixing the survivability and functionality indices blurs the meaning of the objective function. Note that the two indices have clearly different physical interpretations, which are not consistent. Therefore the mixture of the two indices may not give an explicit meaning.
\end{remark}




\section{Formulation of Resilience Control}\label{sec:formulation}

In this section, converters are modeled from the view of DC network side while the topology of DC-SPSs is formulated by extending the spanning tree constraint. Then the two-phase resilience control is presented in detail.

\subsection{Notations}

Index the buses by $1,\cdots,n$ and abbreviate line $\{i,j\}\in\mathcal{E}$ as $i\sim j$. Denote ($i\sim j$ \& $i<j$) as  $i\to j$. 
For each bus $i\in \mathcal{N}$, denote $V_i$ as its voltage. For each $V_i$, let $\underline{V}_i$ and $\overline{V}_i$ denote its lower bound and upper bound, respectively.
$s_i$ indicates a binary parameter, which stands for the available state of load with respect to the functionality of the SPS. Specifically, $s_{i}=1$ if the load is unfaulted; otherwise, $s_{i}=0$.
For each line $i\sim j $, $y_{ij}$ denotes its admittance and $\overline{I}_{ij}$ denotes the upper bound on current. $s_{ij}$ is a binary parameter, which stands for the available state of line $i\sim j$ with respect to the functionality of the SPS. Specifically, $s_{ij}=1$ if the line is unfaulted; otherwise, $s_{ij}=0$.
A letter without subscripts denotes a vector of the corresponding quantities, e.g., $V=[V_i]_{i\in \mathcal{N}}$.

\subsection{Converter Model}
\begin{figure}[htb]
\centering
\includegraphics[width=2.1in]{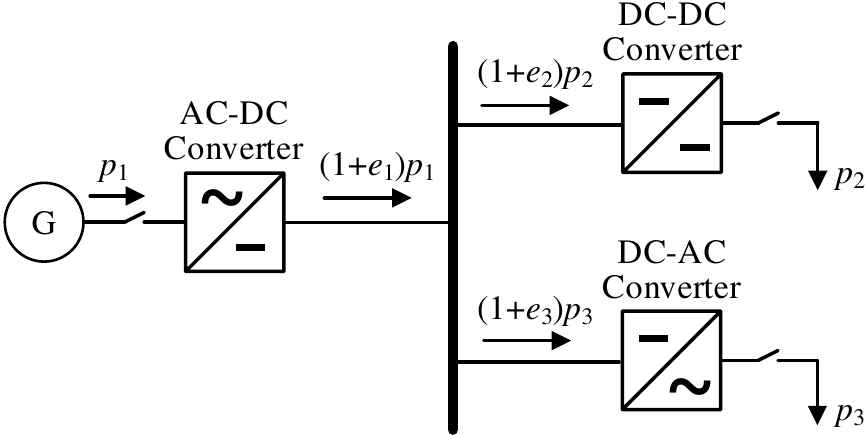}
\caption{Structure of power converters connected to DC-SPS.}
\label{fig:converter}
\end{figure}

In a DC-SPS, AC generators are connected to  DC buses via AC-DC converters. Viewed from the DC network side, the AC-DC converters can be approximately formulated as controllable voltage sources. Since there  are a number of converters in a DC-SPS, the power losses of converters may not be ignored. Denoting $e_i$ as the loss rate of converter $i$, then $(1+e_i)p_i$ represents the power injecting into the DC bus (see Fig. \ref{fig:converter}). Here,  $e_i$ satisfies $0<1+e_i<1$ for AC-DC converters connected to generators, and $1+e_i>1$ for DC-AC and DC-DC converters connected to loads. 

\subsection{Topology Switching Constraints} \label{sec:ring}

\begin{figure}[htb]
	\centering
	\includegraphics[width=3in]{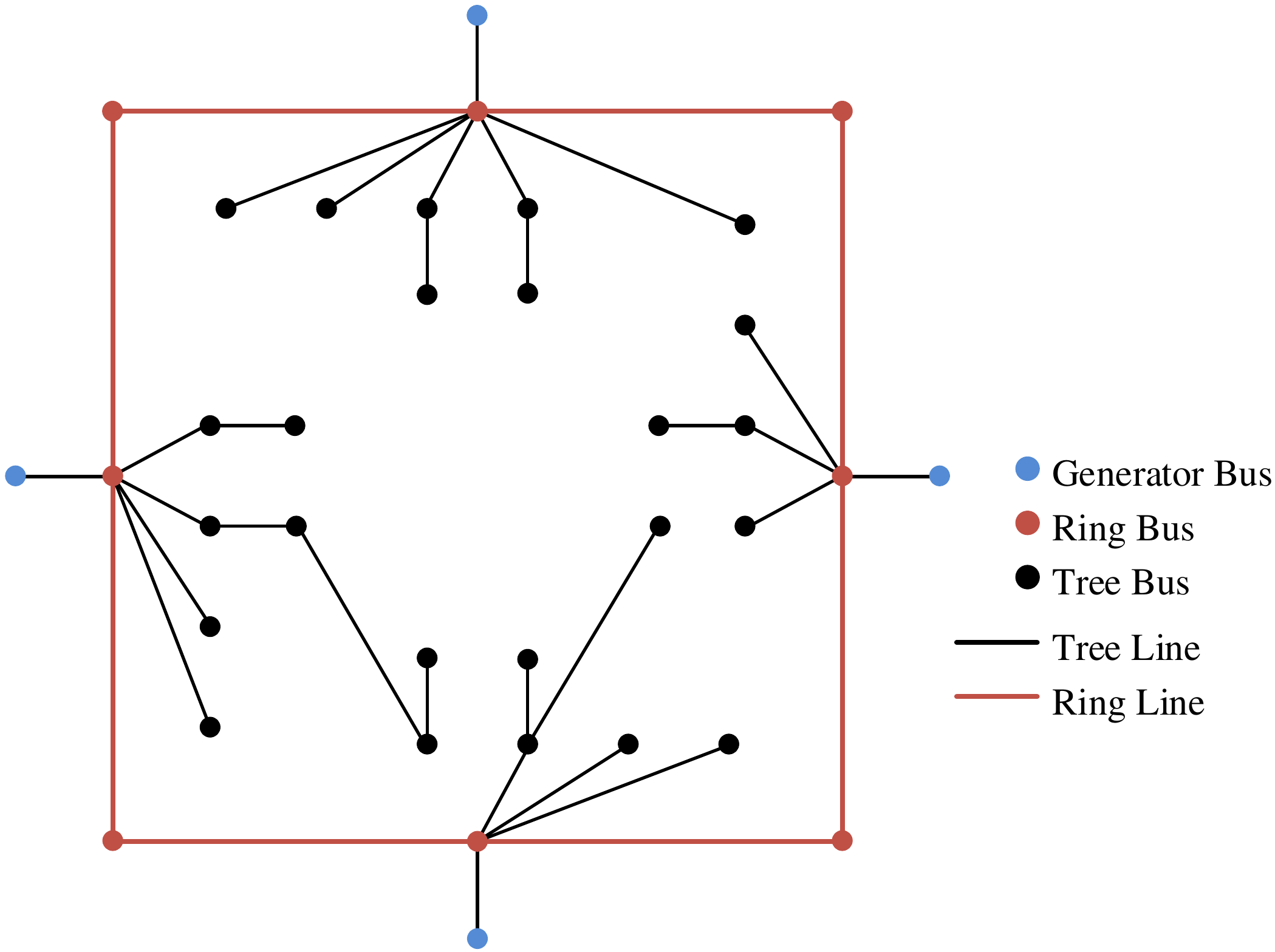}
	\caption{Typical topology of a DC-SPS.}
	\label{fig:ring_topo}
\end{figure}

A typical ring topology of DC-SPSs is shown in Fig \ref{fig:ring_topo}. In this architecture, generators are connected to a ring network through generator switchboards and the subnetworks below the ring are radial. The network flow method \cite{Butler2001} is used to ensure the topologies of the subnetworks connected to the ring network are radial. However, all the nodes in the ring network are merged into one node. It implies that the actions of bus-tie breakers, i.e., ring lines, are not modeled specifically. Thus, the flexibility of network may not be fully exploited. 

Reference \cite{Jabr2012} proposes \lj{the} spanning tree constraint to formulate radial network regardless of the direction of power flow. However, to tackle our problem, it has to be extended to: 1) allow load shedding; 2) introduce a special type of bus, namely ring bus, so that the switching of supply path among multiple generators is allowed.
In this regard, we propose an improved version of the constraint.  

In order to model the ring topology, the set of bus $\mathcal{N}$ is divided into three disjoint subsets: the subset of generator bus $\mathcal{N}^G$, the subset of ring bus $\mathcal{N}^R$ and the subset of tree bus $\mathcal{N}^T$. The generator buses are connected to the ring buses and have no parents. The ring buses are connected to the generator buses or other ring buses. The tree buses are connected in spanning trees.
Similarly, the set of lines $\mathcal{L}$ is divided into two disjoint subsets: the subset of ring lines $\mathcal{L}^R$ and the subset of tree lines $\mathcal{L}^T$. The topology constraints are formulated as below. These constraints model both ring and radial topology, since a radial topology can be viewed as a special type of ring topology with the ring lines disconnected.
\begin{subequations} \label{eq:topo}
\begin{align}
  & \eta_i \le s_i, ~ i\in \mathcal{N}^T \label{eq:eta_s} \\
  & \beta_{ij}\in \{0,1\}, ~ i\sim j \label{eq:beta_binary} \\  
  & \beta_{ij} + \beta_{ji} = \alpha_{ij}, ~ i\sim j \label{eq:beta_alpha} \\
  & 0 \le \alpha_{ij} \le s_{ij}, ~ i\sim j \in \mathcal{L}^T \label{eq:alpha} \\
  & \sum_{j:j\sim i\in \mathcal{L}^T} \beta_{ij} \le 1, ~ i\in \mathcal{N}^T \label{eq:beta} \\
  & \sum_{j:j\sim i\in \mathcal{L}^T} \beta_{ij} \ge \eta_i, ~ i\in \mathcal{N}^T \label{eq:beta_eta} \\
  & \beta_{ij} = 0, ~ i\sim j, i\in \mathcal{N}^G \label{eq:beta_gen} \\ 
  & \beta_{ij} = 0, ~ i\sim j, i\in \mathcal{N}^R, j\in \mathcal{N}^T \label{eq:beta_ring}
\end{align}
\end{subequations}

The state of load is constrained by the availability of load \eqref{eq:eta_s}.
Two binary variables, i.e., $\beta_{ij}$ and $\beta_{ji}$ \eqref{eq:beta_binary}, are introduced for each line $i\sim j$, whose connection state is given by the variable $\alpha_{ij}$. 
Constraint \eqref{eq:beta_alpha} indicates that line $i\sim j$ is connected ($\alpha_{ij}=1$) when either bus $j$ is the parent of bus $i$ ($\beta_{ij}=1$), or bus $i$ is the parent of bus $j$ ($\beta_{ji}=1$). 
It implies that $\alpha_{ij}$ can be treated as a continuous variable constrained by its availability $s_{ij}$ as given in \eqref{eq:alpha}. This formulation empirically accelerates the computation in our case.
The radial condition for the subnetworks below the ring is based on the characteristic of a spanning tree that every node except the root node (generator bus) has only one parent \cite{Jabr2012}. In this context, every tree bus has either one parent connected by a tree line or no parent when it is isolated \eqref{eq:beta}. Additionally, the tree bus always has one parent when the load is switched on \eqref{eq:beta_eta}. The generator buses have no parents \eqref{eq:beta_gen}, while the ring buses only have generator buses or other ring buses as their parents \eqref{eq:beta_ring}.

\subsection{Optimization Formulation of Two-phase Resilience Control}
\subsubsection{Phase I}
The first-phase problem is to maximize the survivability index. Thus the optimization problem is formulated as below.
\begin{subequations}\label{eq:s1}
\begin{align}
  \text{S1:}\ \max ~ & h_1(\eta,p) = \Omega_{sur} - \xi\sum_{i\in \mathcal{N}}(1+e_i)p_i \label{eq:obj1} \\
  \text{over}\ & \alpha, \beta, \eta, \gamma, p, V \nonumber \\
  \text{s.t.}\ & \text{\eqref{eq:eta_s}--\eqref{eq:beta_ring}} \nonumber \\
  & (1+e_i) p_i = \sum_{j:j\sim i}\gamma_{i,ij}\left(\gamma_{i,ij} - \gamma_{j,ij}\right)y_{ij}, ~ i\in \mathcal{N} \label{eq:s1_pi} \\
  & \underline{p}_i s_i \le p_i \le \overline{p}_i s_i, ~ i\in \mathcal{N}^G \label{eq:s1_p_gen} \\
  & \underline{p}_i s_i \eta_i \le p_i \le \overline{p}_i s_i \eta_i, ~ i\in \mathcal{N}^{DV} \label{eq:s1_p_load_variable} \\
  & p_i = \underline{p}_i \eta_i s_i, ~ i\in \mathcal{N}^{DF} \label{eq:s1_p_load_fixed} \\
  & \underline{V}_i \le V_i \le \overline{V}_i, ~ i\in \mathcal{N} \label{eq:s1_v} \\
  & \gamma_{i,ij} = \begin{cases}
  0 & \alpha_{ij} = 0, ~ i\in \mathcal{N}, i\sim j\\
  V_{i} & \alpha_{ij} = 1, ~ i\in \mathcal{N}, i\sim j
  \end{cases} \label{eq:s1_gamma} \\
  & y_{ij}^2(\gamma_{i,ij} - \gamma_{j,ij})^2 \le \overline{I}_{ij}^2 \alpha_{ij}, ~ i\in \mathcal{N}, i\to j \label{eq:s1_current}
\end{align}
\end{subequations}

The total network loss is added in the objective function \eqref{eq:obj1} to make it strictly increasing in $p_i$ for $i\in\mathcal{N}$, and $\xi$ is a small scaling coefficient used to control the impact of power loss on the survivability index $\Omega_{sur}$. This modification is used for proving additional results on exactness later in this paper.
Equality \eqref{eq:s1_pi} denotes the power flow in a DC-SPS with auxiliary variables $\gamma_{i,ij}$ and $\gamma_{j,ij}$. Note that this constraint is non-convex. These auxiliary variables are set as 0 when line $i\sim j$ is disconnected ($\alpha_{ij}=0$), otherwise being equal to $V_i$ and $V_j$ respectively, according to \eqref{eq:s1_gamma}.
The power injection of generator $i$ is constrained by \eqref{eq:s1_p_gen} with the lower bound $\underline{p}_i=0$, since the generators can be turned off.  
There are two types of loads in the model: variable loads \eqref{eq:s1_p_load_variable} and fixed loads \eqref{eq:s1_p_load_fixed}. $\mathcal{N}^{DV}$ denotes the set of variable load buses and $\mathcal{N}^{DF}$ indicates the set of fixed ones. Obviously, one has $\mathcal{N}^D=\mathcal{N}^{DV} \cup \mathcal{N}^{DF}$ and $\mathcal{N}^{DV}  \cap \mathcal{N}^{DF} = \emptyset$.  The lower bound of load $i$ satisfy $\underline{p}_i<0$. Note that in \eqref{eq:s1_p_gen}-\eqref{eq:s1_p_load_fixed}, all $\underline{p}_i$ satisfies $\underline{p}_i\le0$, and the power injections are constrained by their availability.
The nodal voltage $V_i$ is constrained by \eqref{eq:s1_v} with the lower bound $\underline{V}_i>0$. The current is constrained within the upper bound \eqref{eq:s1_current}.

\subsubsection{Phase II}

In the second-phase problem, the variable $\eta$, which stands for the status of load switches, is fixed with the solution obtained in Phase I. To maximize the functionality index, $\Omega_{fun}$, an optimization problem is formulated as below.
\begin{subequations}\label{eq:step2}
\begin{align}
  \text{S2:}\ \max ~ & \Omega_{fun} \nonumber \\
  \text{over}\ & \alpha, \beta, \gamma, p, V \nonumber \\
  \text{s.t.}\ & \text{\eqref{eq:eta_s}--\eqref{eq:beta_ring}}, \text{\eqref{eq:s1_pi}--\eqref{eq:s1_current}} \nonumber
\end{align}
\end{subequations}

The problems S1 and S2 are non-convex due to constraint \eqref{eq:s1_pi}. 
They can  be convexified using SOCP relaxation. We will prove the convexification is exact when the line constraint \eqref{eq:s1_current} is inactive, under some mild assumptions (Assumptions 1 and 2 in Section \ref{sec:convex_relax}).

\begin{remark}
As indicated in \cite{Li}, the exactness of the convex relaxation of OPF in DC networks is independent of network topologies. Later in Section IV, we will further show that this result holds true for \lj{the} SPSs considered in this paper. 
It implies that, to consider all possible topology realizations, one only needs to remove the spanning tree constraints in this section that formulate the ring topology. However, in implementation, a real challenging issue is how to formulate a specific topology into the OPF problem. In this paper, we specifically consider the ring topology due to two reasons: 1) it is a typical topology of DC-SPSs in practice, which is worthy of attention; 2) it is a good example to show how our work can be applied to a network topology with specific requirements.
\end{remark}

\section{Solution Approach}\label{sec:solution}

The proposed framework consists of two mixed-integer non-convex problems, i.e., S1 and S2, where the global optimality may not be guaranteed. In this section, a convex relaxation approach is proposed to convert the problems into MISOCP problems.
For brevity, we only consider S1, since the constraints in S2 are  similar to S1. 

First, by introducing lifted variables, the problem is transformed into an equivalent mixed-integer non-convex counterpart. Second, the non-convex constraints are relaxed, yielding an MISOCP problem.

\subsection{Equivalent Transformation}

First we introduce the lifted variables which satisfy
\begin{subequations}\label{eq:slack_var}
\begin{align}
  \delta_{i,ij} &= \gamma_{i,ij}^2, i\in \mathcal{N}, i\sim j; \label{eq:slack_delta} \\
  W_{ij} &= \gamma_{i,ij}\gamma_{j,ij}, i\sim j; \label{eq:slack_w} \\
  v_i &= V_i^2, i\in\mathcal{N}. \label{eq:slack_v}
\end{align}
\end{subequations}
and define a matrix
\begin{equation*}
R_{ij}:=\begin{bmatrix}
\delta_{i,ij} & W_{ij} \\
W_{ji} & \delta_{j,ij}
\end{bmatrix}
\end{equation*}
for every $i\to j$.
Then the following lemma is given.
\begin{lemma}\label{lemma:trans}
Given the following conditions:
\begin{itemize}
  \item $v_i>0, i\in \mathcal{N}$;
  \item $W_{ij}\ge0, i\to j$;
  \item $W_{ij}=W_{ji}, i\to j$;
  \item $\delta_{i,ij} = \begin{cases}
  0 & \alpha_{ij} = 0, i\in \mathcal{N}, i\sim j \\
  v_i & \alpha_{ij} = 1, i\in \mathcal{N}, i\sim j;
  \end{cases}$
\end{itemize} 
if $\text{rank}(R_{ij})=1$ for all $i\to j$,
then there exists a unique pair of $(V,\gamma)$ which satisfies \eqref{eq:slack_var} and
\begin{subequations} \label{eq:slack_original}
\begin{align}
V_i &> 0,i\in \mathcal{N}; \label{eq:slack_v_0} \\
\gamma_{i,ij} &= \begin{cases}
0 & \alpha_{ij} = 0, i\in \mathcal{N}, i\sim j \\
V_i & \alpha_{ij} = 1, i\in \mathcal{N}, i\sim j. \label{eq:slack_gamma}
\end{cases}
\end{align}
\end{subequations}
Further, such a $(V,\gamma)$ is given by
\begin{subequations}
\begin{align}
V_i &= \sqrt{v_i}, i\in \mathcal{N}; \\
\gamma_{i,ij} &= \sqrt{\delta_{i,ij}}, i\in \mathcal{N}, i\sim j.
\end{align}
\end{subequations}
\end{lemma}

\begin{IEEEproof}
\emph{Existence:} Let $V_i=\sqrt{v_i}$ for $i\in \mathcal{N}$ and $\gamma_{i,ij}=\sqrt{\delta_{i,ij}}$ for $i\in \mathcal{N}$, $i\sim j$. 
It suffices to show that $(V,\gamma)$ satisfies \eqref{eq:slack_w} and \eqref{eq:slack_gamma}.

Since $R_{ij}$ is not full rank, one has
\begin{equation*}
\delta_{i,ij}\delta_{j,ij} - W_{ij}W_{ji} = 0, \quad i\to j.
\end{equation*}
Invoking $W_{ij}\ge0$, we have 
\begin{equation*}
W_{ij} = \sqrt{W_{ij}^2} = \sqrt{W_{ij}W_{ji}} = \sqrt{\delta_{i,ij}\delta_{j,ij}} = \gamma_{i,ij}\gamma_{j,ij}
\end{equation*}
for all $i\sim j$, i.e., $(V,\gamma)$ satisfies \eqref{eq:slack_w}. 
Additionally,
\begin{equation*}
  \gamma_{i,ij} = \sqrt{\delta_{i,ij}} = \begin{cases}
  0 & \alpha_{ij} = 0 \\
  \sqrt{v_i} = V_i & \alpha_{ij} = 1.
  \end{cases}
\end{equation*}
Thus, \eqref{eq:slack_gamma} is satisfied. This completes the proof of existence.

\emph{Uniqueness:} Let $(\tilde{V},\tilde{\gamma})$ denote an arbitrary solution to \eqref{eq:slack_var} and \eqref{eq:slack_original}. It suffices to show that $\tilde{V}_i=\sqrt{v_i}$ for $i\in \mathcal{N}$ and $\tilde{\gamma}_{i,ij}=\sqrt{\delta_{i,ij}}$ for $i\in \mathcal{N}$, $i\sim j$.

Assume $\tilde{V}_i\ne\sqrt{v_i}$ for some $i\in \mathcal{N}$ for the sake of contradiction. Then it follows from \eqref{eq:slack_v} that $\tilde{V}_i=-\sqrt{v_i}<0$, which contradicts with \eqref{eq:slack_v_0}. Thus, $\tilde{V}_i=\sqrt{v_i}$ for all $i\in \mathcal{N}$. 

Similarly, assume $\tilde{\gamma}_{i,ij}\ne\sqrt{\delta_{i,ij}}$ for some $i\in \mathcal{N}$, $i\sim j$. Since it has been proved that $\tilde{V}_i=\sqrt{v_i}$ for $i\in \mathcal{N}$, it follows from \eqref{eq:slack_delta} that 
\begin{equation*}
  \tilde{\gamma}_{i,ij} = -\sqrt{\delta_{i,ij}} = \begin{cases}
  0 & \alpha_{ij} = 0 \\
  -\sqrt{v_i} = -\tilde{V}_i & \alpha_{ij} = 1.
  \end{cases}
\end{equation*}
which contradicts with \eqref{eq:slack_gamma}. Hence, $\tilde{\gamma}_{i,ij}=\sqrt{\delta_{i,ij}}$ for all $i\in \mathcal{N}$ and $i\sim j$. 
This completes the proof of uniqueness.
\end{IEEEproof}

According to Lemma \ref{lemma:trans}, S1 can be transformed equivalently into the following problem.
\begin{subequations}\label{eq:se1}
\begin{align}
  \text{SE1:}\ \max ~ & h_1(\eta,p) \nonumber \\
  \text{over}\ & \alpha, \beta, \eta, \delta, p, v, W \nonumber \\
  \text{s.t.}\ & \text{\eqref{eq:eta_s}--\eqref{eq:beta_ring}}, \text{\eqref{eq:s1_p_gen}--\eqref{eq:s1_p_load_fixed}} \nonumber \\
  & (1+e_i) p_i = \sum_{j:j\sim i}\left(\delta_{i,ij} - W_{ij}\right)y_{ij}, ~ i\in \mathcal{N} \label{eq:se1_p_i} \\
  & \underline{V}_i^2 \le v_i \le \overline{V}_i^2, ~  i\in \mathcal{N} \label{eq:se1_v} \\
  & W_{ij} \ge 0, ~  i\to j \\
  & W_{ij} = W_{ji}, ~ i\to j \\
  & 0 \le \delta_{i,ij} \le \overline{V}_i^2\alpha_{ij}, ~ i\in \mathcal{N}, i\sim j \label{eq:se1_delta_i} \\
  & 0 \le \delta_{j,ij} \le \overline{V}_j^2\alpha_{ij}, ~ j\in \mathcal{N}, i\sim j \\
  & 0 \le v_i - \delta_{i,ij} \le \overline{V}_i^2(1-\alpha_{ij}), ~ i\in \mathcal{N}, i\sim j \\
  & 0 \le v_j - \delta_{j,ij} \le \overline{V}_j^2(1-\alpha_{ij}), ~ j\in \mathcal{N}, i\sim j \label{eq:se1_vj_delta_j} \\
  & R_{ij} \succeq 0, ~ i\to j \label{eq:se1_psd} \\
  & \text{rank}(R_{ij})=1, ~ i\to j  \label{eq:se1_rank} \\
  & y_{ij}^2(\delta_{i,ij} - W_{ij} - W_{ji} + \delta_{j,ij}) \le \overline{I}_{ij}^2 \alpha_{ij}, ~ i\to j \label{eq:se1_current}
\end{align}
\end{subequations}
where the non-convex constraint \eqref{eq:s1_pi} in S1 is transformed to the non-convex constraint \eqref{eq:se1_rank} in SE1. 
$R_{ij}$ is positive semidefinite as shown in \eqref{eq:se1_psd}.
By enforcing constraints \eqref{eq:se1_delta_i}--\eqref{eq:se1_vj_delta_j}, $\delta_{i,ij}$ and $\delta_{j,ij}$ are set as 0 when line $i\sim j$ is disconnected ($\alpha_{ij}=0$), while they take the value of $v_i$ and $v_j$, respectively, when the line is connected ($\alpha_{ij}=1$).

\subsection{Convex Relaxation}\label{sec:convex_relax}

By removing the rank constraint \eqref{eq:se1_rank}, the  non-convex SE1 is transformed to a convex MISOCP problem (named as SR1):
\begin{align}
  \text{SR1:}\ \max ~ & h_1(\eta,p) \nonumber \\
  \text{over}\ & \alpha, \beta, \eta, \delta, p, v, W \nonumber \\
  \text{s.t.}\ & \text{\eqref{eq:eta_s}--\eqref{eq:beta_ring}, \eqref{eq:s1_p_gen}--\eqref{eq:s1_p_load_fixed}, \eqref{eq:se1_p_i}--\eqref{eq:se1_psd}}, \eqref{eq:se1_current}. \nonumber
\end{align}

SR1 is \emph{exact}, provided that its every optimal solution satisfies the rank constraint \eqref{eq:se1_rank}. 
In order to ensure the exactness of convex relaxation when the line constraint \eqref{eq:se1_current} is inactive, throughout the rest of this paper, we make the following assumptions \cite{Li}.
\begin{assumption}\label{assumption1}
  $\overline{V}_1=\overline{V}_2=\cdots\overline{V}_n > 0$.
\end{assumption}

\begin{assumption}\label{assumption2}
 $\sum_{i\in\mathcal{N}}p_i>0$.
\end{assumption}

Assumption \ref{assumption1} requires all the voltages have the same upper bounds, which is reasonable in DC-SPSs. Assumption \ref{assumption2} is trivial as it means the total network loss is positive.

When the line constraint \eqref{eq:se1_current} is inactive, SR1 has the following property.  
\begin{property}\label{property_exact}
    When the line constraint \eqref{eq:se1_current} is inactive, SR1 is exact if Assumptions \ref{assumption1} and \ref{assumption2} hold.
\end{property}

\begin{IEEEproof}
Given predefined parameter $s_{ij}$ and the values of variables $\alpha$, $\beta$, $\eta$, which satisfy \eqref{eq:topo}, the value of $\delta$ is determined according to \eqref{eq:se1_delta_i}--\eqref{eq:se1_vj_delta_j}. Since $\delta_{i,ij}=\delta_{j,ij}=0$ when $\alpha_{ij}=0$, the remaining constraints of SR1 are associated with all connected lines ($\alpha_{ij}=1$) where $\delta_{i,ij}=v_i$ and $\delta_{j,ij}=v_j$. Thus, $R_{ij}$ is transformed into
\begin{equation*}
R'_{ij}:=\begin{bmatrix}
v_i & W_{ij} \\
W_{ji} & v_j
\end{bmatrix}.
\end{equation*}
Accordingly, constraint \eqref{eq:se1_rank} is transformed into 
\begin{equation}
\text{rank}(R'_{ij})=1, ~ i\to j.  \label{eq:sr2_rank}
\end{equation}
Additionally, for any given predefined parameter $s_i$ and binary variable $\eta$ that satisfies \eqref{eq:s1_p_load_fixed}, $p_i$ is constant for $i\in\mathcal{N}^{DF}$. Therefore, 
SR1 is converted into the following problem.
\begin{subequations}
\begin{align}
  \text{SR1':}\ \max ~ & h_1(p) \nonumber \\
  \text{over}\ & p,v,W \nonumber \\
  & p_i = \frac{1}{1 + e_i}\sum_{j:j\sim i}\left(v_i - W_{ij}\right)y_{ij}, i\in \mathcal{N} \label{eq:sr2_p_i} \\
  & \underline{p}_i \le p_i \le \overline{p}_i, i\in \mathcal{N}^G \label{eq:sr2_p_gen} \\
  & \underline{V}_i^2 \le v_i \le \overline{V}_i^2, i\in \mathcal{N} \label{eq:sr2_v} \\
  & W_{ij} \ge 0, i\to j \label{eq:sr2_w_0} \\
  & W_{ij} = W_{ji}, i\to j \label{eq:sr2_w} \\
  & R'_{ij} \succeq 0, ~ i\to j \label{eq:sr2_psd}
\end{align}
\end{subequations}
According to Theorem 2 in \cite{Li}, SR1' is exact. It means that for each predefined parameters $s_{ij}$, $s_i$ and feasible solution $(\alpha,\beta,\eta,\delta)$ to SR1, SR1' is exact. This completes the proof.
\end{IEEEproof}
When the line constraint \eqref{eq:se1_current} is active, the situation turns to be much more complicated, since the rank constraint may not be satisfied. In this context, a slack variable method \cite{Li} can be applied to find a suboptimal solution when the rank constraint is not satisfied. 
After solving SR1, if the solution violates the rank constraint \eqref{eq:se1_rank} for any line $(s\sim t)\in\mathcal{E}$, and any lower bound of $p_s$ and $p_t$ is binding, say, $p_s=\underline{p}_s$, then a corresponding slack variable $\varepsilon_s$ ($\varepsilon_s>0$) is added into \eqref{eq:s1_p_gen} or \eqref{eq:s1_p_load_variable} to reformulate the constraint as
\begin{equation}
  \underline{p}_s \le p_s - \varepsilon_s \le \overline{p}_s
\end{equation}
so that $p_s$ will not reach its lower bound in the next iteration to solve SR1. Additionally, in order to minimize $\varepsilon_s$, it is also added into the objective function as
\begin{equation*}
\max ~ h_1(\eta,p) - \sum_{i\in\hat{\mathcal{N}}} \varepsilon_i
\end{equation*}
where $\hat{\mathcal{N}}$ is the set of buses where the rank constraint \eqref{eq:se1_rank} are violated while the power injection lower bound is binding. 








\subsection{Improving Numerical Stability}

In an optimal solution to SR1, $\delta_{i,ij}$ and $W_{ij}$ may be numerically close to each other, since the range of nodal voltage is small (usually 0.95$\sim$1.05 p.u.) and $\text{Rank}(R_{ij})=1$ is satisfied (implying $v_iv_j=W_{ij}W_{ji}$).
Thus, SR1 is ill-conditioned as \eqref{eq:se1_p_i} needs the subtractions of $\delta_{i,ij}$ and $W_{ij}$, which may cause numerical instability. Such subtractions, however, can be avoided by using a branch flow formulation \cite{Li}, so that the numerical stability is improved. To this end, we define $z_{ij}:=1/y_{ij}$ and adopt alternative variables $P$ and $l$. Then SR1 can be converted into a branch flow model via the mapping $g:(v,\delta,W)\mapsto (v',\delta',P,l)$ defined by:
\begin{equation}\label{eq:convert}
  g:=
\begin{cases}
  v'_i = v_i, & i\in\mathcal{N} \\
  \delta'_{i,ij} = \delta_{i,ij}, & i\in\mathcal{N}, i\to j \\
  P_{ij} = \left(\delta_{i,ij} - W_{ij}\right)y_{ij}, & i\to j \\
  l_{ij} = y_{ij}^2\left(\delta_{i,ij} - W_{ij} - W_{ji} + \delta_{j,ij}\right), & i\to j
\end{cases}
\end{equation}
where $P_{ij}$ denotes the power flow through line $i\to j$, and $l_{ij}$ denotes the magnitude square of the current through line $i\to j$.
With the mapping $g$, SR1 is converted into the following optimization problem with a branch flow formulation.
\begin{subequations}\label{eq:sb1}
\begin{align}
  \text{SB1:} ~ \max ~ & h_1(\eta,p) \nonumber \\
  \text{over} ~ & \alpha, \beta, \eta, \delta, p,P,v,l \nonumber \\
  \text{s.t.} ~ & \text{\eqref{eq:eta_s}--\eqref{eq:beta_ring}, \eqref{eq:s1_p_gen}--\eqref{eq:s1_p_load_fixed}, \eqref{eq:se1_delta_i}--\eqref{eq:se1_vj_delta_j}} \nonumber \\
  & (1 + e_i) p_i = \sum_{j:j\sim i}P_{ij}, i\in \mathcal{N} \label{eq:sb1_p_i} \\
  & \underline{V}_i^2 \le v_i \le \overline{V}_i^2, ~ i\in  \mathcal{N} \label{eq:sb1_v} \\
  & P_{ij} + P_{ji} = z_{ij}l_{ij}, i\to j\\
  & \delta_{i,ij} - \delta_{j,ij} = z_{ij}\left(P_{ij} - P_{ji}\right), i\to j \label{eq:sb1_v_p} \\
  & l_{ij} \ge \frac{P_{ij}^2}{\delta_{i,ij}}, i\sim j \label{eq:sb1_rank} \\
  & l_{ij} \le \overline{I}_{ij}^2 \alpha_{ij}, ~ i\to j \label{eq:sb1_current}
\end{align}
\end{subequations}

SB1  can be solved efficiently by using commercial solvers, such as MOSEK. 
Following the same lines, the Phase-II problem S2 can also be transformed into an MISOCP problem, which is omitted here due to the limit of space.

\section{Case Studies}\label{sec:case}

\begin{figure}[htb]
	\centering
	\includegraphics[width=3.4in]{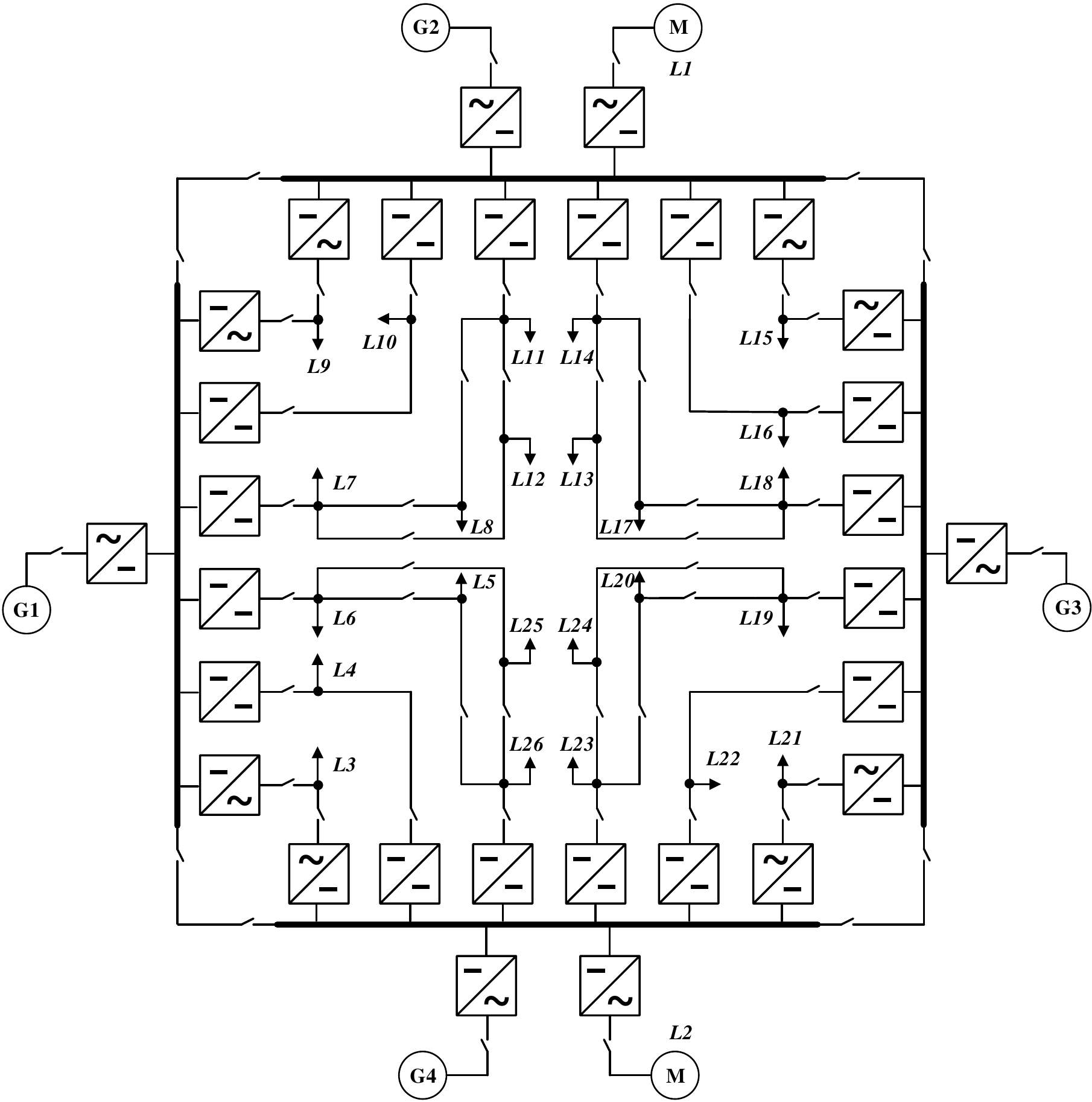}
	\caption{The  DC-SPS for test.}
	\label{fig_ship}
\end{figure}

\begin{figure}[htb]
\centering
\includegraphics[width=2.8in]{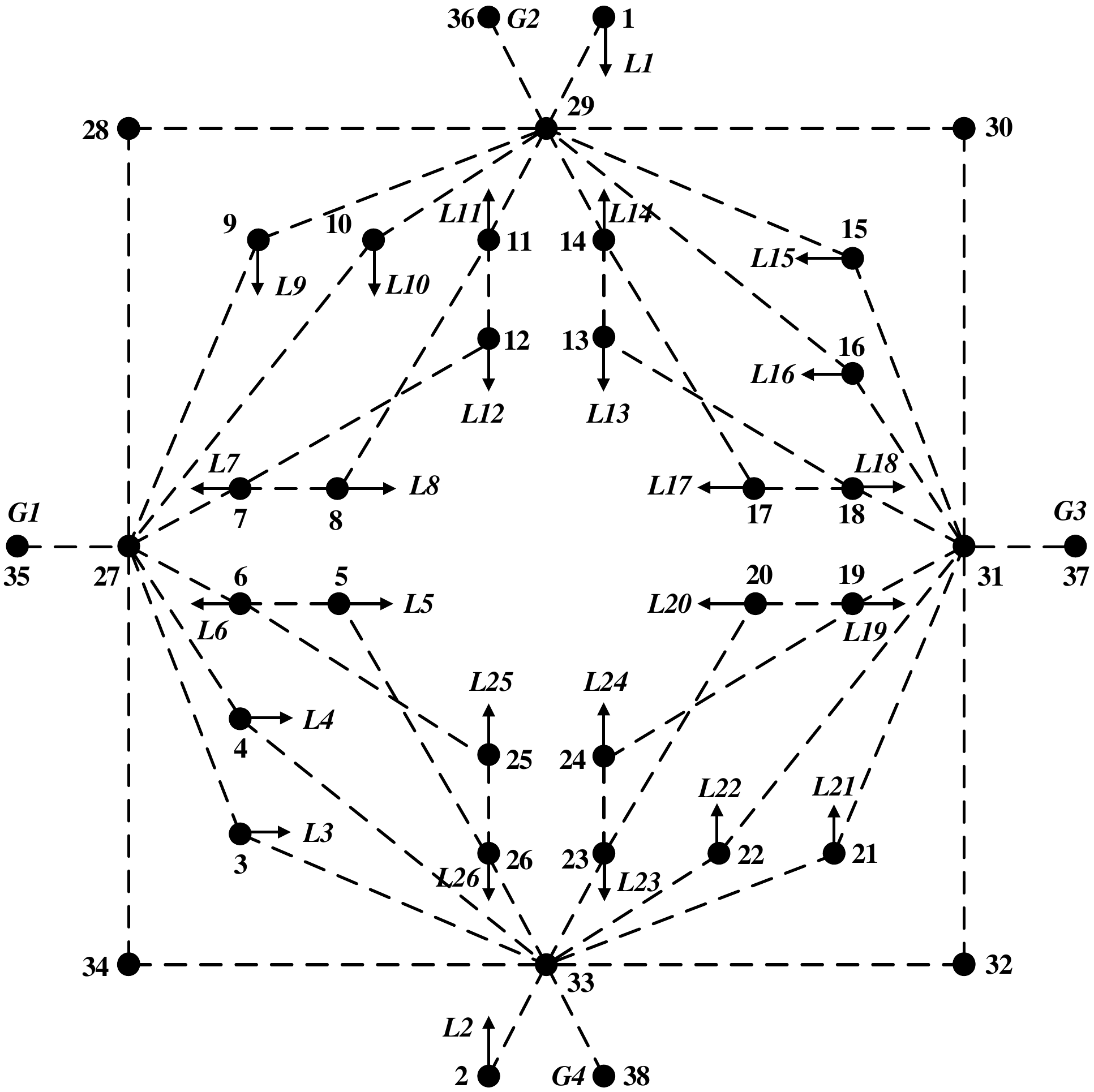}
\caption{Topology of the  DC-SPS for test.}
\label{fig_ship_topo}
\end{figure}

As shown in Fig.\ref{fig_ship}, a DC-SPS is used for testing the proposed methodology. The system contains four generators connected to a ring network. The capacities of G2 and G4 are 4.5 p.u., while the capacities of G1 and G3 are 2 p.u.. The generators are connected to the network via AC-DC converters, while the loads are connected via DC-DC or DC-AC converters. The loss rate $e_i$ is set as -2\% for all the AC-DC converters, and it is set as 2\% for all the DC-DC and DC-AC converters. 
The voltage range is $[0.95, 1.05]$p.u.. The topology of this system is presented in Fig. \ref{fig_ship_topo}, where each dotted line indicates a path. The system consists of 38 buses and 55 paths. Regarding the survivability index, the loads are classified into four priority levels, as listed in Table \ref{tab:priority}. In the literature, there are various ways to select weighting factors\cite{Bose2012,Das2013,Khushalani2008,Bose2010,Amba2009,Mitra2011,Butler-Purry2004}. 
However, there is not a practical standard. In our model, the weighting factors are related to the number of buses in an SPS. Without loss of generality, we assume the weighting factor of the lowest priority level (level 4) is 1. Other weighting factors satisfy \eqref{eq:weight}.
Other parameters can be found in Appendix.

The MISOCP problem is solved on a laptop with two 2.60 GHz Intel Core i5 processors and 8 GB of RAM. The convergence tolerance is set as $10^{-7}$. The branch flow formulation \eqref{eq:sb1} is adopted for all the tests.

\begin{table}[htbp]
  \centering
  \caption{Priority Levels of Loads}
    \begin{tabular}{ccc}
    \hline
    Priority Level & $\kappa_i^k$ & Load Bus Number \\
    \hline
    1     & 729   & 1, 2 \\
    2     & 81    & 3, 5, 9, 12, 15, 17, 21, 24 \\
    3     & 9     & 4, 8, 10, 13, 16, 20, 22, 25 \\
    4     & 1     & 6, 7, 11, 14, 18, 19, 23, 26 \\
    \hline
    \end{tabular}%
  \label{tab:priority}%
\end{table}%

\subsection{Survivability and Functionality}

\begin{figure}[htb]
\centering
\includegraphics[width=2.8in]{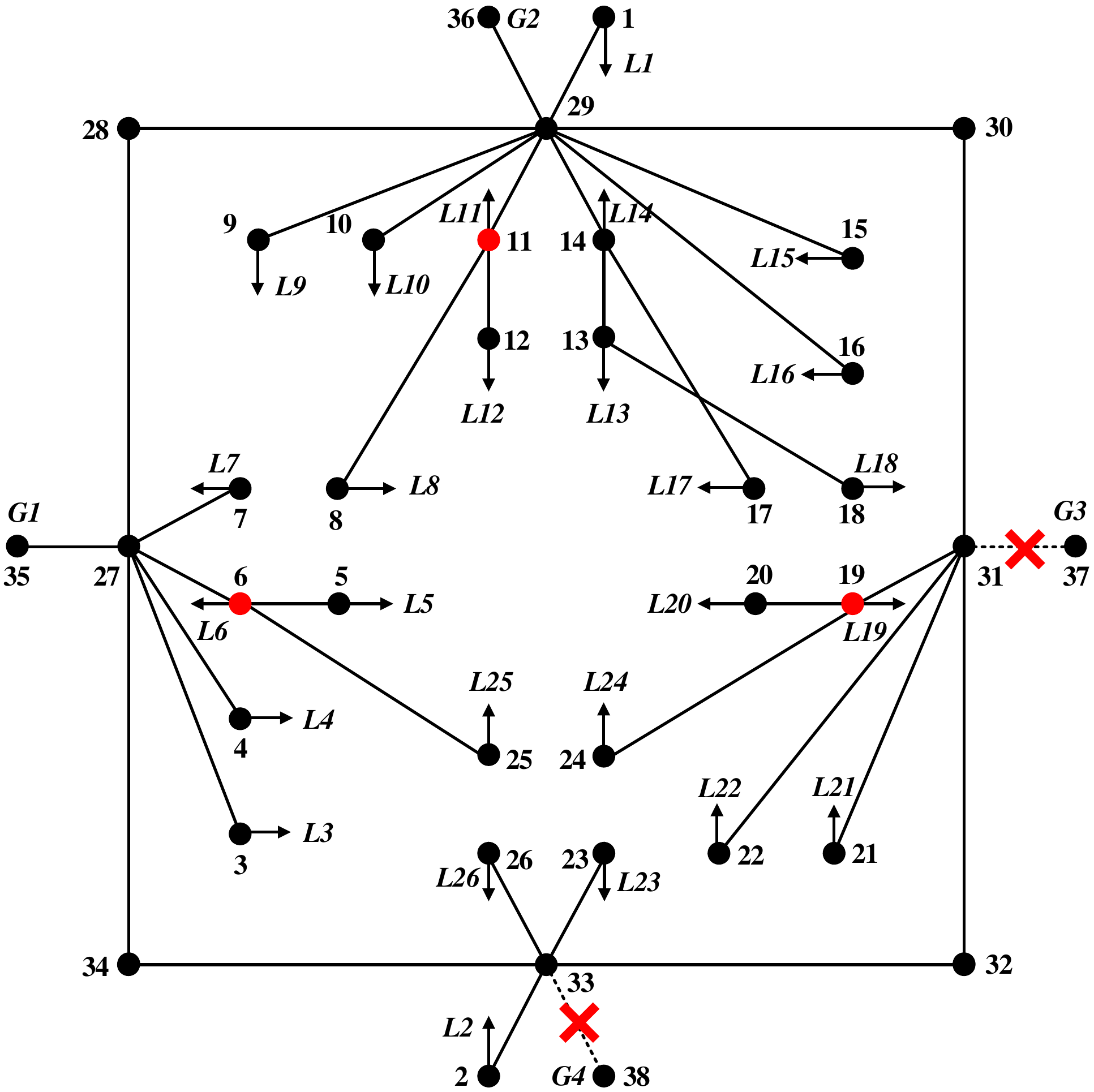}
\caption{Topology of the DC-SPS after G3 and G4 are disconnected.}
\label{fig_survive}
\end{figure}

As the generators directly affect the survivability of DC-SPS, we first test different faults on the lines connected to generators. The results are listed in Table \ref{tab:generator}. 
When G1 or G2 is disconnected, the survivability index $\Omega_{sur}$ still achieves its maximum value, i.e., 1. It means that no load has to be switched off. 
However, due to the loss of generator capacity, the functionality index $\Omega_{fun}$ is less than 1, implying that some of the loads may not be fully supplied. How much these loads are supplied is determined by solving the optimization problem of maximizing the total functionality. Since G2 has a larger capacity than G1, $\Omega_{fun}$ is smaller when G2 is disconnected. 
When G3 and G4 are disconnected simultaneously, more loads have to be switched off. The proposed method ensures the loads with the lowest priority (see Table \ref{tab:priority}) are switched off first. Additionally, the topology after G3 and G4 are disconnected is depicted in Fig \ref{fig_survive}, where the red crosses indicate the positions of faults. After faults, the loads supplied by G3 and G4 are supplied by G1 and G2 via the ring lines.

\begin{table}[htbp]
  \centering
  \caption{Results of the Faults on the Lines Connected to Generators}
    \begin{tabular}{ccccc}
    \hline
    Fault Lines & Disconnected  & $\Omega_{sur}$ & $\Omega_{fun}$ & Bus No. of \\
    & Generator & & & Switch-off Loads \\
    \hline
    27-35        & G1      & 1  & 0.991 & - \\
    29-36        & G2      & 1  & 0.883 & - \\
    31-37, 33-38 & G3, G4  & 0.999  & 0.680 & 6, 11, 19 \\
    \hline
    \end{tabular}%
  \label{tab:generator}%
\end{table}%

\subsection{Re-sectionalization}

\begin{figure}[htb]
\centering
\includegraphics[width=2.8in]{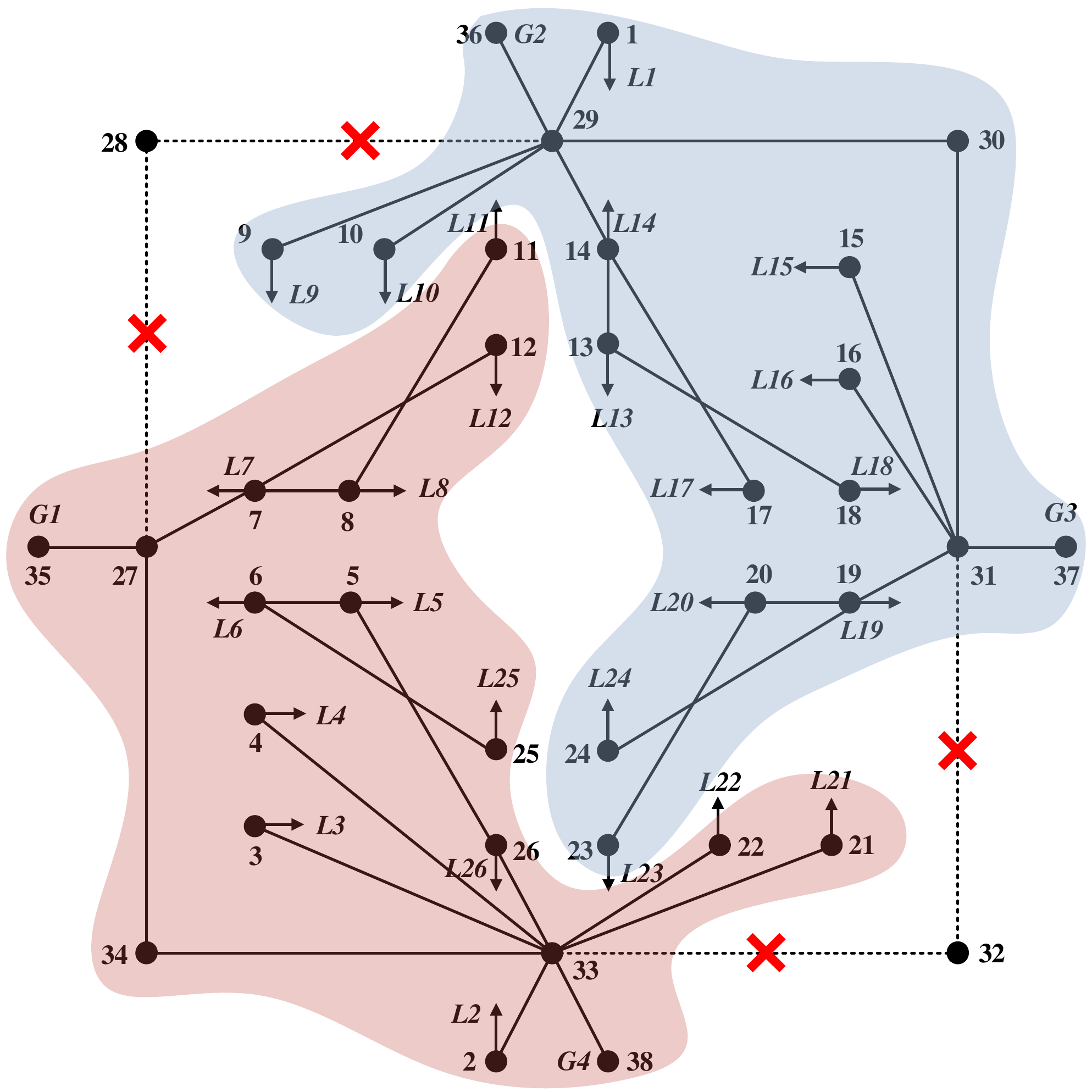}
\caption{Topology of the DC-SPS after the faults on the ring lines.}
\label{fig_sec1}
\end{figure}

\begin{figure}[htb]
\centering
\includegraphics[width=2.8in]{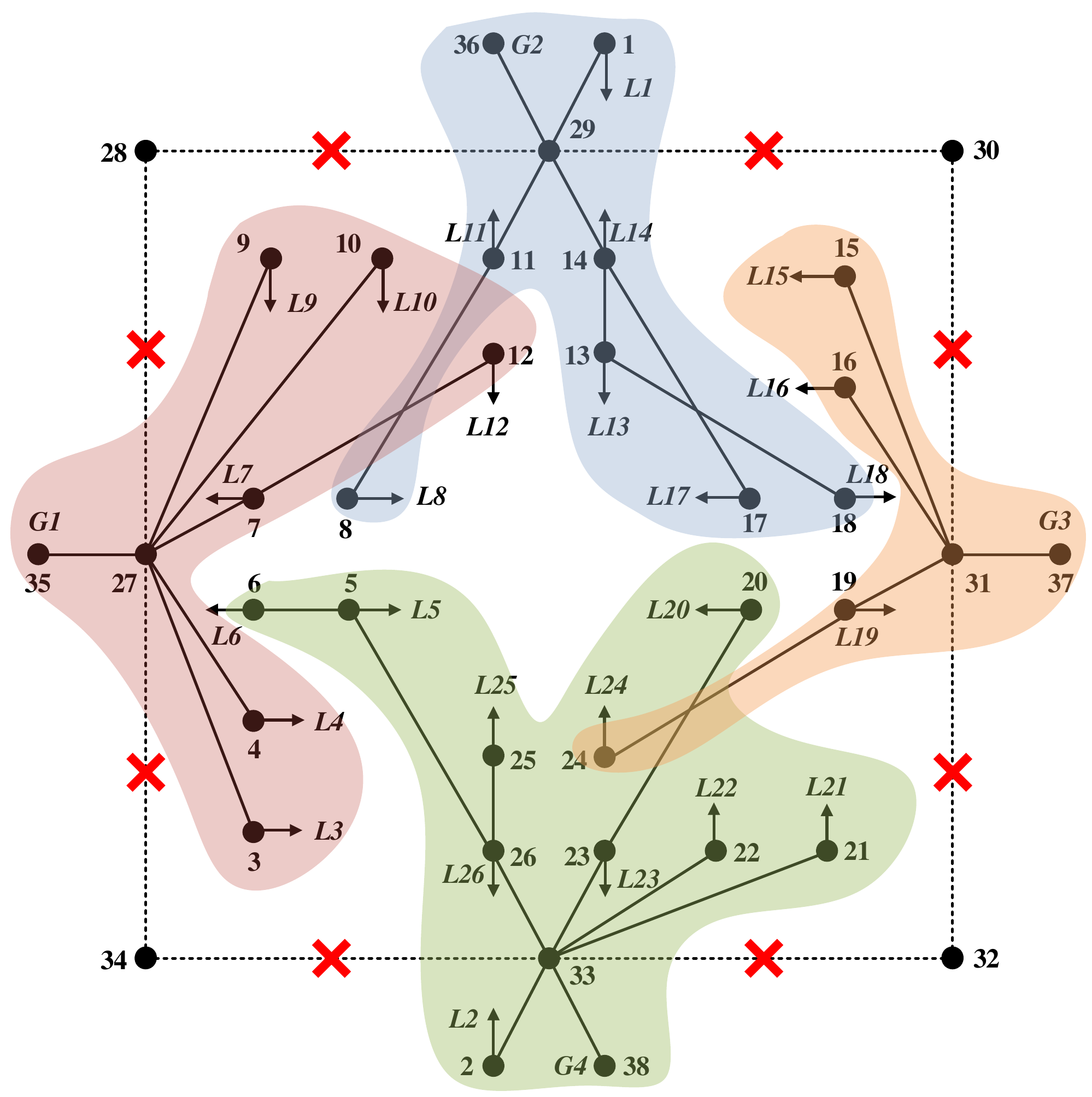}
\caption{Topology of the DC-SPS after the faults on the ring lines.}
\label{fig_sec3}
\end{figure}

In order to exploit the flexibility of re-sectionalization, different faults on the ring lines are simulated. The resulting topologies are shown in Figs. \ref{fig_sec1}$\sim$\ref{fig_sec3}. It is observed that the whole system is re-sectionalized and the lines remain radial topologies below the ring. In all these cases, all the loads can be fully satisfied.

\begin{figure}[htb]
\centering
\includegraphics[width=2.8in]{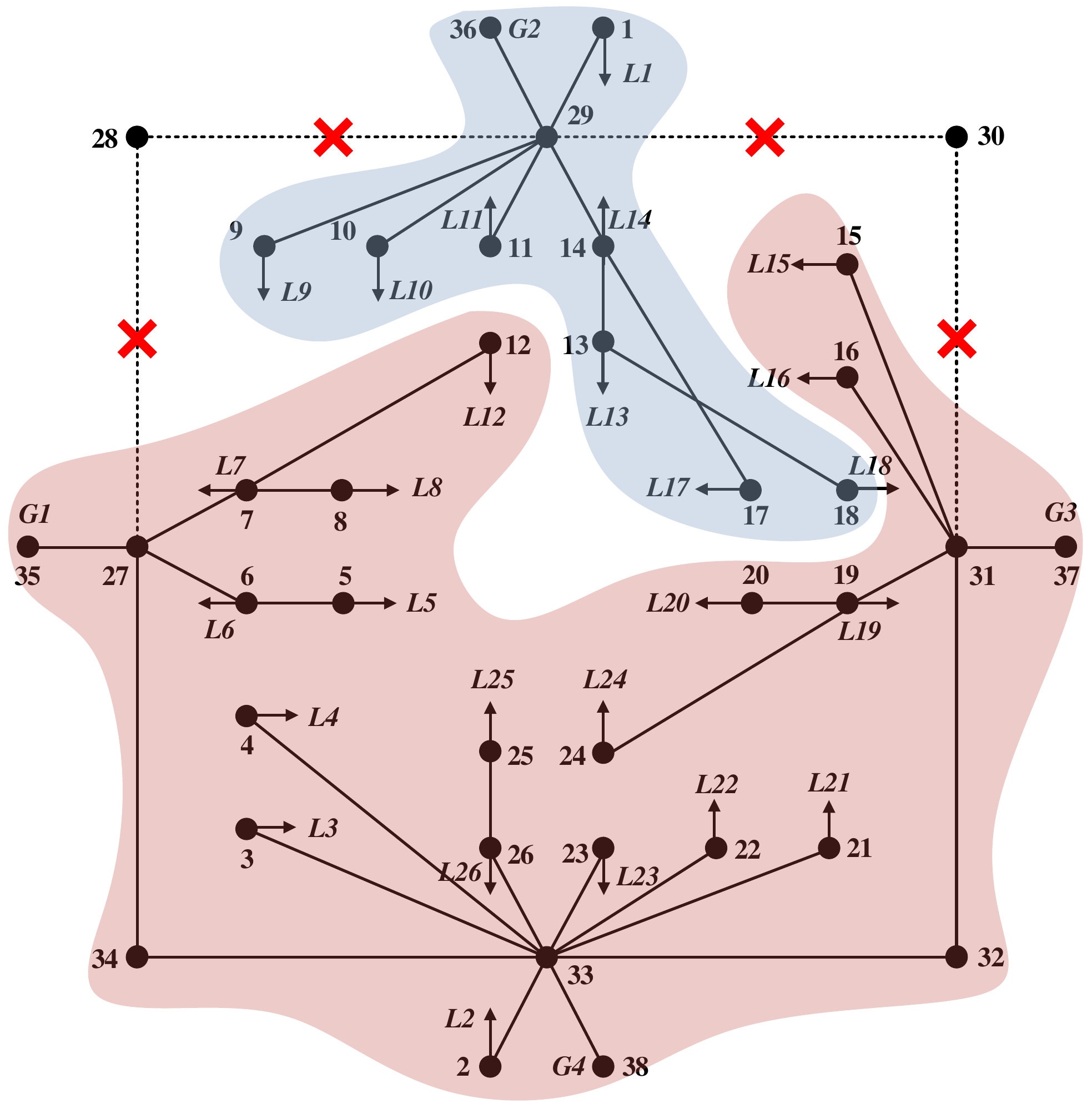}
\caption{Topology of the DC-SPS after the faults on the ring lines.}
\label{fig_sec2}
\end{figure}

\subsection{Reconfiguration}

\begin{figure}[htb]
\centering
\includegraphics[width=2.8in]{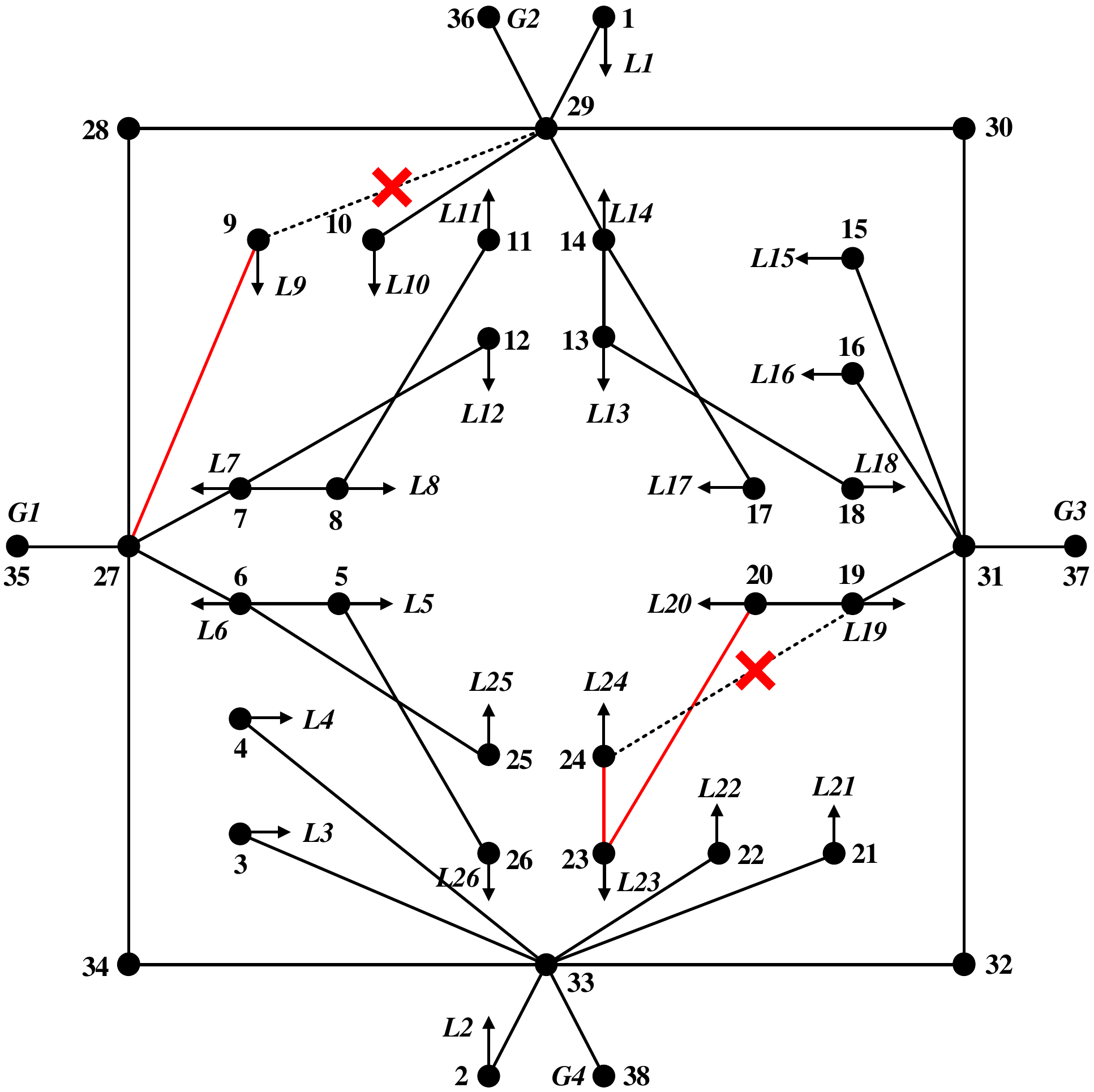}
\caption{Topology of the DC-SPS after the faults on the tree lines.}
\label{fig_switch}
\end{figure}

Figure \ref{fig_switch} shows the topology after the faults occur upon two lines. Before the fault, load 9 is satisfied via line 9-29. After the fault, the supply path is switched to line 9-27. Load 24 is supplied via line 19-24 before the fault, and the supply path is switched to lines 23-24 and 20-23 after restoration.
In this case, the converter loss accounts for about 95\% of the total power loss. It indicates that the converter loss should not be ignored in DC-SPSs when a large amount of load is supplied via converters.

\subsection{Optimality and Computational Efficiency}

Seven cases are simulated to test the computational efficiency, where the faults occur upon different positions, including lines connected generator (GL), lines connected to load (LL), and ring lines (RL). 
\lj{Three} approaches are compared:

\begin{description}
	\item[\textbf{A1:}] The proposed MISOCP approach solved using MOSEK 7.0.0.75;
	\item[\textbf{A2:}] The original non-convex problem (S1) solved using DICOPT and IPOPT 3.11;
  \item[\textbf{A3:}] The original non-convex problem (S1) solved using DICOPT and CONOPT 3.15L.
\end{description}


\begin{table*}[htbp]
  \centering
  \scriptsize
  \caption{Results of Fault Combinations}
    \begin{tabular}{ccccccccccccc}
    \hline
    Case  &       &       & \multicolumn{3}{c}{A1} & \multicolumn{3}{c}{A2} & \multicolumn{3}{c}{A3} \\
    No.   & Fault Combinations & Fault Lines & $\Omega_{sur}$ & $\Omega_{fun}$ & Time (s) & $\Omega_{sur}$ & $\Omega_{fun}$ & Time (s) & $\Omega_{sur}$ & $\Omega_{fun}$ & Time (s) \\
    \hline
    1     & GL    & 33-38 & 1     & 0.883 & 1.42  & 1     & 0.883 & 13.00 & 0.958 & 0.976 & 3.89 \\
    2     & RL    & 27-28, 28-29, 30-31 & 1     & 1     & 2.54  & 1     & 1     & 15.44 & -     & -     & - \\
    3     & LL    & 5-6, 14-29, 19-20 & 1     & 1     & 1.96  & 1     & 1     & 32.14 & 1     & 1     & 3.14 \\
    4     & GL+RL & 27-28, 27-34, 27-35 & 1     & 0.991 & 10.46 & 1     & 0.991 & 34.58 & 0.537 & 1     & 4.80 \\
    5     & GL+LL & 3-27, 3-33, 7-8, 33-38 & 0.963 & 0.896 & 1.68  & 0.963 & 0.896 & 52.90 & 0.911 & 1     & 2.65 \\
    6     & RL+LL & 13-14, 28-29, 29-30 & 1     & 1     & 2.18  & 1     & 1     & 18.70 & 0.949 & 0.975 & 2.94 \\
    7     & GL+RL+LL & 5-26, 27-35, 29-30 & 1     & 0.991 & 3.29  & 1     & 0.991 & 29.39 & -     & -     & - \\
    \hline
    \end{tabular}%
  \label{tab:efficiency}%
\end{table*}%

As shown in Table \ref{tab:efficiency}, whereas the computation time of A3 is comparable with A1 in most of  cases, A3 fails to obtain the optimal solution in Phase I except case 3. Much worse, A3 stops on maximum iteration cycle in case 2, and fails to converge in case 7 due to the inherent non-convexity of the problem. When A2 is adopted, it achieves the same optimal solutions as A1 in all cases. However, A1 is much faster than A2. The acceleration can be as much as an order of magnitude. 
These cases empirically demonstrate that the proposed method can find the optimal solutions efficiently.

\section{Concluding Remarks}\label{sec:conclusion}

In this paper, we have proposed the survivability index and functionality index, rendering a two-phase resilience control methodology for DC-SPSs. Benefiting from the proposed topology formulation, the flexibility of network can be fully exploited. Various topologies, including radial, ring, zonal, can be optimally utilized when different faults occur. Converter loss is also taken into account as it contributes a significant proportion to the total power loss, particularly when \lj{most of the} load\lj{s} \lj{are} supplied.

The proposed model is inherently a mixed-integer non-convex programming problem. 
To meet the practical requirement of fast resilience control in an emergency, an MISOCP-based algorithm is presented to solve the problem efficiently.
Numerical tests empirically verify the effectiveness and efficiency of the proposed algorithm.

Although this paper is focused on resilience control of DC-SPSs, we would like to highlight that DC-SPS is merely a typical example of stand-alone DC networks. Since \lj{either} a grid-connected DC microgrid or a DC distribution network, \lj{which} contains a substation bus with an unconstrained power injection and a fixed nodal voltage, can be regarded as a special case of \lj{a} stand-alone DC network. In this sense, the proposed resilience control methodology can be readily extended to various kinds of DC microgrids and DC distribution networks.
 
It is worth mentioning that, in this study, we have not considered stability constraints, nor the sequence of control actions. These issues are also crucial for practical applications, which are our future works.

\ifCLASSOPTIONcaptionsoff
  \newpage
\fi

\bibliographystyle{IEEEtran}
\bibliography{library}

\begin{thebibliography}{10}
\providecommand{\url}[1]{#1}
\csname url@samestyle\endcsname
\providecommand{\newblock}{\relax}
\providecommand{\bibinfo}[2]{#2}
\providecommand{\BIBentrySTDinterwordspacing}{\spaceskip=0pt\relax}
\providecommand{\BIBentryALTinterwordstretchfactor}{4}
\providecommand{\BIBentryALTinterwordspacing}{\spaceskip=\fontdimen2\font plus
\BIBentryALTinterwordstretchfactor\fontdimen3\font minus
  \fontdimen4\font\relax}
\providecommand{\BIBforeignlanguage}[2]{{%
\expandafter\ifx\csname l@#1\endcsname\relax
\typeout{** WARNING: IEEEtran.bst: No hyphenation pattern has been}%
\typeout{** loaded for the language `#1'. Using the pattern for}%
\typeout{** the default language instead.}%
\else
\language=\csname l@#1\endcsname
\fi
#2}}
\providecommand{\BIBdecl}{\relax}
\BIBdecl

\bibitem{Jin2016a}
Z.~Jin, G.~Sulligoi, R.~Cuzner, L.~Meng, J.~C. Vasquez, and J.~M. Guerrero,
  ``{Next-generation shipboard DC power system: introduction smart grid and dc
  microgrid technologies into maritime electrical netowrks},'' \emph{IEEE
  Electrif. Mag.}, vol.~4, no.~2, pp. 45--57, 2016.

\bibitem{Reed2012}
G.~F. Reed, B.~M. Grainger, A.~R. Sparacino, and Z.-H. Mao, ``{Ship to grid:
  medium-voltage DC concepts in theory and practice},'' \emph{IEEE Power Energy
  Mag.}, vol.~10, no.~6, pp. 70--79, 2012.

\bibitem{Butler1999}
K.~L. Butler, N.~Sarma, and V.~R. Prasad, ``{A new method of network
  reconfiguration for service restoration in shipboard power systems},''
  \emph{IEEE Transm. Distrib. Conf.}, vol.~2, pp. 658--662, 1999.

\bibitem{Butler2001}
K.~L. Butler, N.~D.~R. Sarma, and V.~{Ragendra Prasad}, ``{Network
  reconfiguration for service restoration in shipboard power distribution
  systems},'' \emph{IEEE Trans. Power Syst.}, vol.~16, no.~4, pp. 653--661,
  2001.

\bibitem{Butler-Purry2004}
K.~L. Butler-Purry and N.~D.~R. Sarma, ``{Self-healing reconfiguration for
  restoration of naval shipboard power systems},'' \emph{IEEE Trans. Power
  Syst.}, vol.~19, no.~2, pp. 754--762, 2004.

\bibitem{Butler-Purry2005}
K.~L. Butler-Purry, ``{Multi-agent technology for self-healing shipboard power
  systems},'' \emph{Proc. 13th Int. Conf. Intell. Syst. Appl. to Power Syst.
  ISAP'05}, vol. 2005, pp. 207--211, 2005.

\bibitem{Solanki2005}
J.~M. Solanki and N.~N. Schulz, ``{Using intelligent multi-agent systems for
  shipboard power systems reconfiguration},'' \emph{Proc. 13th Int. Conf.
  Intell. Syst. Appl. to Power Syst. ISAP'05}, vol. 2005, pp. 212--214, 2005.

\bibitem{Huang2007}
K.~Huang, D.~A. Cartes, and S.~K. Srivastava, ``{A multiagent-based algorithm
  for ring-structured shipboard power system reconfiguration},'' \emph{IEEE
  Trans. Syst. Man Cybern. Part C Appl. Rev.}, vol.~37, no.~5, pp. 1016--1021,
  2007.

\bibitem{Feliachi2006a}
A.~Feliachi, K.~Schoder, S.~Ganesh, and H.-J. Lai, ``{Distributed control
  agents approach to energy management in electric shipboard power system},''
  in \emph{2006 IEEE Power Eng. Soc. Gen. Meet.}, 2006, p. 6 pp.

\bibitem{Jing2009}
H.~Jing, Z.~Xiaofeng, and C.~Yan, ``{Multiobjective optimal service restoration
  for shipboard power system},'' \emph{Proc. - 2009 IEEE Int. Conf. Intell.
  Comput. Intell. Syst. ICIS 2009}, vol.~1, pp. 644--650, 2009.

\bibitem{Shariatzadeh2016}
F.~Shariatzadeh, N.~Kumar, and A.~Srivastava, ``{Optimal control algorithms for
  reconfiguration of shipboard microgrid distribution system using intelligent
  techniques},'' \emph{IEEE Trans. Ind. Appl.}, vol. 9994, no.~c, pp. 1--1,
  2016.

\bibitem{Kumar2007}
N.~Kumar, A.~K. Srivastava, and N.~N. Schulz, ``{Shipboard power system
  restoration using binary particle swarm optimization},'' \emph{39th IEEE
  North Am. Power Symp.}, pp. 164--169, 2007.

\bibitem{Mitra2011}
P.~Mitra and G.~K. Venayagamoorthy, ``{Implementation of an intelligent
  reconfiguration algorithm for an electric ship's power system},'' \emph{IEEE
  Trans. Ind. Appl.}, vol.~47, no.~5, pp. 2292--2300, 2011.

\bibitem{Hari2010}
S.~Hari, K.~Vuppalapati, and A.~K. Srivastava, ``{Application of ant colony
  optimization for reconfiguration of shipboard power system},'' vol.~2, no.~3,
  pp. 119--131, 2010.

\bibitem{Khushalani2008}
S.~Khushalani, J.~Solanki, and N.~Schulz, ``{Optimized restoration of combined
  ac/dc shipboard power systems including distributed generation and islanding
  techniques},'' \emph{Electr. Power Syst. Res.}, vol.~78, no.~9, pp.
  1528--1536, sep 2008.

\bibitem{Amba2009}
T.~Amba, K.~L. Butler-Purry, and M.~Falahi, ``{Genetic algorithm based damage
  control for shipboard power systems},'' in \emph{2009 IEEE Electr. Sh.
  Technol. Symp.}\hskip 1em plus 0.5em minus 0.4em\relax Baltimore, MD, USA:
  IEEE, 2009, pp. 242--252.

\bibitem{Bose2010}
S.~Bose, S.~Pal, C.~Scoglio, B.~Natarajan, S.~Das, and N.~Schulz, ``{Analysis
  of optimized reconfiguration of power system for electric ships},'' in
  \emph{North Am. Power Symp. 2010}.\hskip 1em plus 0.5em minus 0.4em\relax
  Arlington, TX, USA: IEEE, 2010, pp. 1--7.

\bibitem{Bose2012}
S.~Bose, S.~Pal, B.~Natarajan, C.~M. Scoglio, S.~Das, and N.~N. Schulz,
  ``{Analysis of optimal reconfiguration of shipboard power systems},''
  \emph{IEEE Trans. Power Syst.}, vol.~27, no.~1, pp. 189--197, 2012.

\bibitem{Das2013}
S.~Das, S.~Bose, {Siddhartha Pal}, N.~N. Schulz, C.~M. Scoglio, and
  B.~Natarajan, ``{Dynamic reconfiguration of shipboard power using
  reinforcement learning},'' \emph{IEEE Trans. Power Syst.}, vol.~28, no.~2,
  pp. 669--676, 2013.

\bibitem{Hijazi2015}
H.~Hijazi and S.~Thi{\'{e}}baux, ``{Optimal distribution systems
  reconfiguration for radial and meshed grids},'' \emph{Int. J. Electr. Power
  Energy Syst.}, vol.~72, pp. 136--143, 2015.

\bibitem{Gan2014}
L.~Gan and S.~H. Low, ``{Optimal power flow in direct current networks},''
  \emph{IEEE Trans. Power Syst.}, vol.~29, no.~6, pp. 2892--2904, 2014.

\bibitem{Li}
J.~Li, F.~Liu, Z.~Wang, S.~H. Low, and S.~Mei, ``{Optimal power flow in
  stand-alone DC microgrids},'' \emph{IEEE Trans. Power Syst.}, pp. 1--1, 2018.

\bibitem{Nagarajan2016}
H.~Nagarajan, E.~Yamangil, R.~Bent, P.~{Van Hentenryck}, and S.~Backhaus,
  ``{Optimal resilient transmission grid design},'' in \emph{2016 Power Syst.
  Comput. Conf.}\hskip 1em plus 0.5em minus 0.4em\relax Genoa, Italy: IEEE,
  2016, pp. 1--7.

\bibitem{Barnes2017}
\BIBentryALTinterwordspacing
A.~Barnes, H.~Nagarajan, E.~Yamangil, R.~Bent, and S.~Backhaus, ``{Tools for
  improving resilience of electric distribution systems with networked
  microgrids},'' pp. 1--8, 2017. [Online]. Available:
  \url{http://arxiv.org/abs/1705.08229}
\BIBentrySTDinterwordspacing

\bibitem{Jabr2012}
R.~A. Jabr, R.~Singh, and B.~C. Pal, ``{Minimum loss network reconfiguration
  using mixed-integer convex programming},'' \emph{IEEE Trans. Power Syst.},
  vol.~27, no.~2, pp. 1106--1115, 2012.

\end{thebibliography}

\appendix
\begin{table}[htbp]
  \centering
  \tabcolsep2pt
  \scriptsize
  \begin{threeparttable}
  \caption{Load and Generator Parameters}
    \begin{tabular}{cccccc}
    \hline
    Load/Generator & $\underline{p}_i$ (p.u.) & $\overline{p}_i$ (p.u.) & Load/Generator & Min (p.u.) & Max (p.u.) \\
    \hline
    L1    & -1.000 & -0.500 & L16   & -0.225 & -0.113 \\
    L2    & -1.000 & -0.500 & L17   & -0.770 & - \\
    L3    & -0.168 & -     & L18   & -0.220 & -0.110 \\
    L4    & -0.208 & -0.104 & L19   & -0.750 & -0.375 \\
    L5    & -0.670 & -     & L20   & -0.860 & -0.430 \\
    L6    & -0.610 & -0.305 & L21   & -0.135 & - \\
    L7    & -0.520 & -0.260 & L22   & -0.217 & -0.109 \\
    L8    & -0.630 & -0.315 & L23   & -0.027 & -0.014 \\
    L9    & -0.163 & -     & L24   & -0.690 & - \\
    L10   & -0.163 & -0.082 & L25   & -0.230 & -0.115 \\
    L11   & -0.620 & -0.310 & L26   & -0.031 & -0.016 \\
    L12   & -0.300 & -     & G1    & 0     & 2 \\
    L13   & -0.810 & -0.405 & G2    & 0     & 4.5 \\
    L14   & -0.210 & -0.105 & G3    & 0     & 2 \\
    L15   & -0.123 & -     & G4    & 0     & 4.5 \\
    \hline
    \end{tabular}%
  \label{tab:load_gen}
      \begin{tablenotes}
        \item [a] The type of fixed load only has a lower bound.
    \end{tablenotes}
    \end{threeparttable}   
\end{table}%

\begin{table}[htbp]
  \centering
  \scriptsize
  \tabcolsep2pt
  \caption{Line Parameters}
    \begin{tabular}{cccccc}
    \hline
    From Bus & To Bus & Resistance (p.u.)  & From Bus & To Bus & Resistance (p.u.)  \\
    \hline
    1     & 29    & 2.1*1e-4 & 16    & 31    & 2.2*1e-4 \\
    2     & 33    & 1.2*1e-4 & 17    & 18    & 2.6*1e-4 \\
    3     & 27    & 3.3*1e-4 & 18    & 31    & 1.1*1e-4 \\
    3     & 33    & 2.1*1e-4 & 19    & 20    & 3.3*1e-4 \\
    4     & 27    & 4.3*1e-4 & 19    & 24    & 5.0*1e-4 \\
    4     & 33    & 1.1*1e-4 & 19    & 31    & 2.2*1e-4 \\
    5     & 6     & 5.2*1e-4 & 20    & 23    & 1.2*1e-4 \\
    5     & 26    & 1.3*1e-4 & 21    & 31    & 4.1*1e-4 \\
    6     & 25    & 2.1*1e-4 & 21    & 33    & 2.5*1e-4 \\
    6     & 27    & 3.5*1e-4 & 22    & 31    & 5.2*1e-4 \\
    7     & 8     & 1.2*1e-4 & 22    & 33    & 1.3*1e-4 \\
    7     & 12    & 4.4*1e-4 & 23    & 24    & 2.3*1e-4 \\
    7     & 27    & 2.1*1e-4 & 23    & 33    & 4.2*1e-4 \\
    8     & 11    & 3.1*1e-4 & 25    & 26    & 1.1*1e-4 \\
    9     & 27    & 5.2*1e-4 & 26    & 33    & 5.3*1e-4 \\
    9     & 29    & 1.3*1e-4 & 27    & 28    & 2.2*1e-4 \\
    10    & 27    & 2.4*1e-4 & 27    & 34    & 2.4*1e-4 \\
    10    & 29    & 4.5*1e-4 & 27    & 35    & 1.5*1e-4 \\
    11    & 12    & 3.1*1e-4 & 28    & 29    & 2.2*1e-4 \\
    11    & 29    & 5.3*1e-4 & 29    & 30    & 2.1*1e-4 \\
    13    & 14    & 1.2*1e-4 & 29    & 36    & 2.3*1e-4 \\
    13    & 18    & 4.4*1e-4 & 30    & 31    & 2.0*1e-4 \\
    14    & 17    & 2.4*1e-4 & 31    & 32    & 2.3*1e-4 \\
    14    & 29    & 1.5*1e-4 & 31    & 37    & 4.1*1e-4 \\
    15    & 29    & 3.1*1e-4 & 32    & 33    & 2.1*1e-4 \\
    15    & 31    & 5.0*1e-4 & 33    & 34    & 2.3*1e-4 \\
    16    & 29    & 6.4*1e-4 & 33    & 38    & 4.1*1e-4 \\
    \hline
    \end{tabular}%
  \label{tab:line_para}%
\end{table}%

\begin{IEEEbiographynophoto}{Jia Li}
received the B.S. degree in electrical engineering from Tsinghua University, Beijing, China, in 2012. He is currently pursuing the Ph.D. degree at the Department of Electrical Engineering, Tsinghua University. His current research interests include power system operation and planning under uncertainty.
\end{IEEEbiographynophoto}

\begin{IEEEbiographynophoto}{Feng Liu}
(M'10) received the B.Sc. and Ph.D. degrees in electrical engineering from Tsinghua University, Beijing, China, in 1999 and 2004, respectively. 
Dr. Liu is currently an Associate Professor of Tsinghua University. From 2015 to 2016, he was a visiting associate at the California Institute of Technology, CA, USA. His research interests include power system stability analysis, optimal control and robust dispatch, game theory and learning theory and their applications to smart grids. He is the author/coauthor of more than 100 peer-reviewed technical papers and two books, and holds more than 20 issued/pending patents. He was a guest editor of IEEE Transactions on Energy Conversion.
\end{IEEEbiographynophoto}

\begin{IEEEbiographynophoto}{Ying Chen}
(M'07) received the B.E. and Ph.D. degrees in electrical engineering from Tsinghua University, Beijing, China, in 2001 and 2006, respectively. He is currently an Associate Professor with the Department of Electrical Engineering and Applied Electronic Technology, Tsinghua University. His research interests include parallel and distributed computing, electromagnetic transient simulation, cyber-physical system modeling, and cyber security of smart grid.
\end{IEEEbiographynophoto}

\begin{IEEEbiographynophoto}{Chengcheng Shao}
(M'17) received the B.S. and Ph.D. degree in electrical engineering from Xi’an Jiaotong University, Xi’an, China, in 2011 and 2017. He is currently an Assistant Professor at Xi’an Jiaotong University. His research interests include the integration of wind power and load dispatch.
\end{IEEEbiographynophoto}

\begin{IEEEbiographynophoto}{Guanqun Wang}
(St. M'12-M'16) received his B.S. and Ph.D. degrees in electrical engineering from Tsinghua University, China, in 2007 and 2012, respectively. He received Ph.D. degree in electrical engineering at Washington State University, USA, in 2016. His research interests include wide area monitoring and control, and power system economic dispatch. He is a staff transmission planning engineer in Burns \& McDonnell.
\end{IEEEbiographynophoto}

\begin{IEEEbiographynophoto}{Yunhe Hou}
(M'08-SM'15) received the B.E. and Ph.D. degrees in electrical engineering from Huazhong University of Science and Technology, Wuhan, China, in 1999 and 2005, respectively. He was a Post-Doctoral Research Fellow at Tsinghua University, Beijing, China, from 2005 to 2007, and a Post-Doctoral Researcher at Iowa State University, Ames, IA, USA, and the University College Dublin, Dublin, Ireland, from 2008 to 2009. He was also a Visiting Scientist at the Laboratory for Information and Decision Systems, Massachusetts Institute of Technology, Cambridge, MA, USA, in 2010. Since 2017, he has been a Guest Professor with Huazhong University of Science and Technology, China. He joined the faculty of the University of Hong Kong, Hong Kong, in 2009, where he is currently an Associate Professor with the Department of Electrical and Electronic Engineering. Dr. Hou is an Editor of the IEEE Transactions on Smart Grid and an Associate Editor Journal of Modern Power Systems and Clean Energy.
\end{IEEEbiographynophoto}

\begin{IEEEbiographynophoto}{Shengwei Mei}
(SM'06-F'14) received the B.Sc. degree in mathematics from Xinjiang University, Urumqi, China, the M.Sc. degree in operations research from Tsinghua University, Beijing, China, and the Ph.D. degree in automatic control from Chinese Academy of Sciences, Beijing, China, in 1984, 1989 and 1996, respectively. He is currently a Professor at Tsinghua University. His research interests include power system analysis and control, game theory, and its application in power systems.
\end{IEEEbiographynophoto}

\end{document}